\documentclass{article}
\usepackage{latexsym,amssymb}

\begin{document}

\begin{center}
{\bf On Endo-trivial Modules for $p$-Solvable Groups}\\
\medskip

\emph{Gabriel Navarro\\
Departament d'\`Algebra\\
Facultat de Matem\`atiques\\
Universitat de Val\`encia\\
46100 Burjassot\\ Val\`encia\\
SPAIN\\
E-mail: gabriel.navarro@uv.es\\
 \& \\
Geoffrey R. Robinson,\\
Institute of Mathematics\\
University of Aberdeen\\
Aberdeen\\
AB24 3UE\\
SCOTLAND\\
E-mail: g.r.robinson@abdn.ac.uk\\
July 2010}
\end{center}

\medskip
\noindent {\bf Introduction:} In this note, we will prove a
conjecture of J. Carlson, N. Mazza and J. Th\'evenaz [1], namely, we will
prove that if $G$ is a finite $p$-nilpotent group which contains a
non-cyclic elementary Abelian $p$-subgroup and  $k$ is an
algebraically closed field of characteristic $p,$ then all simple
endo-trivial $kG$-modules are $1$-dimensional. In fact, we do rather
more: we prove the analogous result directly in the case that $G$ is
$p$-solvable and contains an elementary Abelian $p$-subgroup of order
$p^{2}.$ Carlson, Mazza and Th\'evenaz had reduced the proof of this
result for $p$-solvable $G$ to the $p$-nilpotent case, (and had
proved the result in the solvable case), but our method is somewhat
different. Our proof does require the classification of finite
simple groups. Specifically, we require the well-known fact that the outer automorphism group
of a  finite simple group of order prime to $p$ has cyclic Sylow $p$-subgroups (see, for example,
Theorem 7.1.2 of Gorenstein, Lyons and Solomon, [3]).

\medskip
Let us recall that a $kG$-module $M$ is endo-trivial if $M \otimes
M^{*} \cong k \oplus N,$ where $N$ is a projective $kG$-module. If
$|G|$ is divisible by $p,$ any endo-trivial $kG$-module has
dimension prime to $p.$ The vertex of any indecomposable
endo-trivial $kG$-module is a Sylow $p$-subgroup of $G$. We remark
that if $M$ is an endo-trivial $kG$-module which is not
$1$-dimensional, then a Sylow $p$-subgroup of $G$ acts faithfully on
$M \otimes M^{*}$ and hence acts faithfully on $M.$

\medskip
In the first Lemma, we summarize some
properties of endo-trivial modules which are for the most part
well-known. We recall that a subgroup $H$ of a finite group $G$ is
said to be strongly $p$-embedded if $p$ divides $|H|$ and $p$ does not
divide $|H \cap H^{g}|$ for each $g \in G \backslash H.$

\medskip
\noindent {\bf Lemma 1:} \emph{ Let $M$ be an endo-trivial module
for a finite group $X.$ Then:\\
\noindent i) Every non-projective summand of ${\rm Res}_{Z}^{X}(M)$
is endo-trivial for each subgroup
$Z$ of $X.$\\
\noindent ii) $M$ has a unique non-projective indecomposable summand
which is itself endo-trivial.\\
\noindent iii) If $M \cong U \otimes V$ for $kX$-modules $U$ and $V,$ then both $U$ and $V$ are
endo-trivial.\\
\noindent iv) If $M \cong {\rm Ind}_{Y}^{X}(L)$ for some proper subgroup
$Y$ of $X,$ and some $kY$-module $L,$ then $L$ is endo-trivial and
$Y$ is strongly $p$-embedded in $X.$}

\medskip
\noindent {\bf PROOF:} The first two parts are clear. To prove iii),
notice that $U$ and $V$ each have dimension prime to $p.$ Hence $U
\otimes U^{*} =  k \oplus S$ and $V \otimes V^{*} = k \oplus T,$
where $S$ and $T$ are $kX$-modules. Then $M \otimes M^{*} = k \oplus
S \oplus T \oplus (S \otimes T),$ so that $S$ and $T$ must both be
projective.

\medskip
\noindent iv) Notice that $L$ is isomorphic to a non-projective
direct summand of ${\rm Res}^{X}_{Y}(M)$ in this case, so that $L$
is endo-trivial. Since $M$ has dimension prime to $p,$ we see
that $Y$ must contain a Sylow $p$-subgroup of $X.$

\medskip
Now $$M\otimes M^{*} \cong {\rm Ind}_{Y}^{X}[ L\otimes {\rm
Res}^{X}_{Y}(M^{*})],$$ so that $M \otimes M^{*}$ has a direct
summand isomorphic to ${\rm Ind}^{X}_{Y}(L \otimes L^{*}),$ and in
particular, a direct summand isomorphic to ${\rm Ind}_{Y}^{X}(k).$
Since $M$ is endo-trivial, this implies that the only non-projective
indecomposable summand of ${\rm Ind}_{Y}^{X}(k)$ is $k.$ By the
Mackey formula (applied to the restriction of this permutation
module to $Y$)  this implies that ${\rm Ind}_{Y \cap Y^{x}}^{Y}(k)$
is projective for each $x \in X \backslash Y,$ so that $Y \cap
Y^{x}$ is a $p^{\prime}$-subgroup for each $ x \in X \backslash Y$
and $Y$ is strongly $p$-embedded in $X.$ 

\medskip
\noindent {\bf Remark:} The converse of part iv) is also true: if an endo-trivial
$kY$-module is induced from the strongly
$p$-embedded subgroup $Y$ of $X$, then the resulting $kX$-module
is also endo-trivial. This is almost immediate from the proof of iv) above and Mackey's Theorem.

\medskip
\noindent {\bf Corollary:} \emph{ Let $X$ be a $p$-solvable finite
group containing an elementary Abelian subgroup of order $p^{2}.$
Then no endo-trivial $kX$-module is induced from a proper subgroup
of $X.$}

\medskip
\noindent {\bf Proof:} Let $Q$ be a Sylow $p$-subgroup of $X.$ If
such a module were induced from a proper subgroup $Y$ of $X,$ then
$Y$ would be strongly $p$-embedded in $X$ (and may be chosen to
contain $Q$) by the previous Lemma. Let $Z = O_{p^{\prime},p}(X).$
Then $X = O_{p^{\prime}}(X)N_{X}(Q \cap Z).$ Now, as $Q$ contains an
elementary Abelian subgroup of order $p^{2}$, by 6.2.4 of
Gorenstein, [2], for example, we have $$O_{p^{\prime}}(X) \leq \langle
C_{X}(u): u \in Q^{\#}\rangle \leq Y,$$ and we also have $N_{X}(Q
\cap Z) \leq Y.$ Hence $X \leq Y,$ contrary to the fact that $Y$ is
proper.

\medskip
The following Lemma is well-known, but we include its proof for completeness:

\medskip
\noindent {\bf Lemma 2:} \emph{Let $G = \langle x \rangle N$ be a finite $p$-nilpotent group with Sylow $p$-subgroup $\langle x \rangle$
of order $p$ and normal $p$-complement $N.$ Suppose that $V$ is a simple endo-trivial $kG$-module, and let $W$ be its Green correspondent
for $N_{G}(\langle x \rangle) (= C_{G}(x)).$ Then all indecomposable summands of ${\rm Res}^{C_{G}(x)}_{C_{N}(x)}(W)$
are isomorphic and $1$-dimensional.}

\medskip
\noindent {\bf Proof:} Since $V$ has dimension prime to $p,$ the restriction of $V$ to $N$ is simple.
Notice that $C_{G}(x) = \langle x \rangle \times C_{N}(x),$ and that every indecomposable $kC_{G}(x)$-module
is expressible as a tensor product $A \otimes B$ where $C_{N}(x)$ acts trivially on $A$ and $x$ acts indecomposably
on $A,$ while $x$ acts trivially on $B$ and $C_{N}(x)$ acts irreducibly on $B.$

\medskip
Now $W$ is indecomposable, and is also endo-trivial. Writing $W$ in
the above fashion as $A \otimes B,$ both $A$ and $B$ are
endo-trivial by part iii) of Lemma 1. Then $B$ is $1$-dimensional,
since $\langle x \rangle$ acts trivially on $B.$ (Notice also that
${\rm dim}_{k}(A) \leq p,$ so we either have ${\rm dim}_{k}(A) =
p-1$ or ${\rm dim}_{k}(A) = 1).$

\medskip
We recall that a component of a finite group $X$ is a subnormal quasi-simple
subgroup of $X.$ Distinct components of $X$ centralize each other,
and all components of $X$ centralize the Fitting subgroup $F(X).$
The central product of the components of $X$ is denoted by $E(X).$
The following Lemma is probably well-known, but we include a proof.

\medskip
\noindent {\bf Lemma 3:} \emph{Let $G$ be a perfect finite group
with $G =E(G)$ and with $Z = Z(G)$ a cyclic $p^{\prime}$-group. Suppose that $G \lhd H$
for some finite group $H$ with $Z \leq Z(H).$ Suppose further that the components of $G$ are all conjugate
within $H$ and that the element $x$ of order $p$ in $H$ permutes the components of $G$ semi-regularly by
conjugation. Then $C_{G}(x)$ is isomorphic to a central product of
components of $G,$ one from each $\langle x \rangle$-orbit. In
particular, $C_{G}(x)$ is perfect.}

\medskip
\noindent {\bf Proof:} We first note that $Z$ is contained in each
component of $G.$ For suppose that $L$ is a component of $G$ and $ W
= L \cap Z < Z.$ Then all components of $G/W$ are simple, and $G/W$
has a non-trivial Abelian direct factor $Z/W,$ so is not perfect, a
contradiction. Let $L_{1},L_{2},\ldots, L_{n}$ be representatives
for the $\langle x \rangle$-orbits of components of $G.$ Let $T =
L_{1}L_{2}\ldots L_{n},$ a central product of mutually centralizing
components. Then $T,T^{x},\ldots,T^{x^{p-1}}$ are mutually
centralizing, since no two of them contain a common component. We
may thus define a homomorphism $\phi: T \to C_{G}(x)$ via $t\phi =
tt^{x}\ldots t^{x^{p-1}}.$ In the case that $Z = 1,$ this is clearly
a surjection. When $Z \neq 1,$ we have $Z \leq T\phi$ since $x$ acts
trivially on $Z$ and $Z$ is a $p^{\prime}$-group. Also, we have
(by, for example, 5.3.15 of Gorenstein [2]) $C_{H/Z}(xZ) = C_{H}(x)/Z.$ The analysis in the $Z =1$ applies to $G/Z,$ so
that $C_{G/Z}(xZ)$ is clearly isomorphic to $T/Z$ and hence $C_{G}(x)$ is
isomorphic to $T$ since $T$ injects into $C_{G}(x).$

\medskip
\noindent {\bf Theorem:} \emph{ Let $X$ be a $p$-solvable group which
contains an elementary Abelian subgroup $Q$ of order $p^{2}$
and let $V$ be a simple endo-trivial $kX$-module. Then
$V$ is $1$-dimensional.}

\medskip
\noindent {\bf Proof:} If possible, choose a counterexample $(X,V)$ so that
first ${\dim}_{k}(V),$ then $|X|,$ are minimized. Then $V$ is a faithful $kP$-module, where $P$ is
a Sylow $p$-subgroup of $G.$ But $V$ is simple, so that $O_{p}(X)$ acts trivially on $V.$ Hence $O_{p}(X) = 1.$
More generally, the kernel of the action of $X$ on $V$ is a $p^{\prime}$-group, so that $V$ is
a faithful $kX$-module by minimality. We know that $V$ is a primitive $kX$-module by the Corollary above. This enables us to perform
standard Clifford-theoretic reductions, and the endo-trivial condition turns out to be compatible with these reductions.

\medskip
Let $Y$ be a normal subgroup of $X$ minimal subject to strictly containing $Z(X).$  Since $Y/Z(X)$ is a minimal normal subgroup of $X/Z(X),$ we know that $Y/Z(X)$ is a direct product of simple groups. If $Y/Z(X)$ is Abelian, then $Y^{\prime} \leq Z(X)$ and $Y$ is nilpotent. Notice that $Y/Z(X)$ is not a $p$-group as
$O_{p}(X) =1.$ Let $U$ be an irreducible summand of ${\rm Res}^{X}_{Y}(V).$ Since $V$ is primitive and
$Y$ is non-central, the isomorphism type of $U$ is $X$-stable, but $U$ is not $1$-dimensional. The usual Clifford-theoretic
construction yields a $p^{\prime}$-central extension ${\hat X}$ of $X$ such that $U$ extends to a simple $k{\hat X}$-module,
and such that $V \cong U \otimes W$ as $k{\hat X}$-module (in fact $W \cong {\rm Hom}_{{\hat Y}}(U,V)$).
By part iii) of Lemma 1, both $U$ and $W$ are endo-trivial as $k{\hat X}$-modules, since $V$ is also endotrivial
as $k{\hat X}$-module (for ${\hat X}$ acts as $X$ does on $V \otimes V^{*}).$ The Sylow $p$-subgroups of $X$ and of ${\hat X}$ are clearly isomorphic.
If neither $U$ nor $W$ is one dimensional, we have a contradiction to the minimal choice of $(X,V).$ Hence $W$ must be one dimensional, as $U$ is not.

\medskip
Hence ${\rm dim}_{k}(U) = {\rm dim}_{k}(V),$ so that $V$ restricts irreducibly to $Y.$ Hence
$Z(Y) \leq Z(X)= C_{X}(Y)$ by Schur's Lemma. Now $Y$ is either nilpotent of class $2$ 
 or else is the central product of $Z(X)$ with a single
conjugacy class of components, each of order prime to $p$. By the minimal choice of
$(X,V),$ we now have $X = YQ.$ If $Y$ is not nilpotent of class $2,$ then $Y^{\prime} = E(Y)$ still acts irreducibly on $V,$
so the minimal choice of $(X,V)$ gives $X = E(Y)Q  = E(X)Q$ in that case.

\medskip
Suppose that $Y$ is nilpotent. This case was deal with by Carlson, Mazza and Th\'evenaz in [1], but we provide a different argument
to dispose of it. We know that $Y/Z(Y)$ is a minimal normal subgroup of $X/Z(X).$ For any $a \in Q^{\#},$ we have $Z(Y) \leq C_{Y}(a) \lhd Y.$
Furthermore, $C_{Y}(a)$ is $Q$-invariant, so $C_{Y}(a) \lhd YQ = X.$ However, $C_{Y}(a) \neq Y,$ as $Y$ acts irreducibly on $V$
and $a$ has order $p.$ Thus $C_{Y}(a) = Z(Y).$ Since $a$ was arbitrary, and $Z(Y) \neq Y = \langle C_{Y}(a): a \in Q^{\#}\rangle,$
this is a contradiction. Hence $Y = E(Y).$

\medskip
We have already remarked that $Q$ acts
faithfully on $V.$ Hence no element of $Q$ can centralize $Y,$ as
$Y$ acts irreducibly on $V.$ Since $Y$ is a $p^{\prime}$-group and $Q$ centralizes $Z(Y),$ the action of $Q$ on $Y/Z(Y)$ is faithful.
Hence $Y$ is not quasi-simple, for otherwise the outer automorphism group of $Y/Z(Y)$ has cyclic Sylow $p$-subgroups. Since
$Y/Z(Y)$ is a minimal normal subgroup of $X/Z(X),$ the components of $Y$
are transitively permuted under conjugation by $Q.$ Hence there is
an element $a \in Q$ which acts semi-regularly by conjugation on the
components of $Y.$ Then $C_{Y}(a)$ is perfect by Lemma 3. In
particular, there is no non-trivial $1$-dimensional simple $kC_{Y}(a)$-module.

\medskip
However, by Lemma 2, the Green correspondent of ${\rm
Res}^{X}_{Y\langle a \rangle}(V)$ for $C_{X}(a)$ lies over a
1-dimensional module for $C_{Y}(a),$ so lies over the trivial module
of $C_{Y}(a).$ Hence that Green correspondent lies in the principal
block of $C_{X}(a).$ By Brauer's Third Theorem (and the
compatibility between Green correspondence and Brauer
correspondence), ${\rm Res}^{X}_{Y\langle a \rangle}(V)$ lies in the
principal block of $kY\langle a \rangle,$ a contradiction, as $Y$
acts faithfully on $V$ and $Y = O_{p^{\prime}}(Y\langle a \rangle).$

\medskip
\begin{center}
{\bf REFERENCES}
\end{center}

\medskip
\noindent [1] J. F. Carlson, N. Mazza and J. Th\'evenaz, \emph
{Endotrivial modules for $p$-solvable groups}, to appear in
Transactions of the American Mathematical Society.

\medskip
\noindent [2] D. Gorenstein, \emph{Finite groups}, Second edition,
Chelsea Publishing Co., New York, 1980.

\medskip
\noindent [3] D. Gorenstein, R. Lyons and R. Solomon, \emph{The Classification of the Finite Simple Groups, Number 3},
American Mathematical Society Mathematical Surveys and Monographs, Volume 40, Number 3, AMS, Providence, 1998.

\end{document}